\documentclass [12pt,oneside,a4paper,mathscr]{amsart}

\usepackage {amscd,amsmath,amssymb,euscript}

\newtheorem{thm}{Theorem}[section]
\newtheorem{cor}[thm]{Corollary}
\newtheorem{lemma}[thm]{Lemma}
\newtheorem{prop}[thm]{Proposition}

\theoremstyle{definition}
\newtheorem{defn}[thm]{Definition}
\newtheorem{eg}[thm]{{\bf Example}}

\theoremstyle{remark}
\newtheorem*{pf}{Proof}
\newtheorem*{ack}{Acknowledgements}
\newtheorem*{nota}{Notation}

\newcommand{\tensor}{\otimes}
\newcommand{\HH}{H}
\newcommand{\dual}{\vee}
\newcommand{\T}{\operatorname{T}}
\newcommand{\A}{\mathcal A}
\newcommand{\B}{\mathcal B}
\newcommand{\oom}{\omega}
\newcommand{\om}{\Omega}
\newcommand{\eps}{\epsilon}
\newcommand{\Hom}{\operatorname{Hom}}
\newcommand{\lHom}{\operatorname{{\mathcal H}om}}
\newcommand{\lra}{\longrightarrow}
\newcommand{\lRa}[1]{\stackrel{#1}{\lra}}
\newcommand{\lla}{\longleftarrow}
\newcommand{\lLa}[1]{\stackrel{\,#1}{\lla}}

\newcommand{\isom}{\cong}
\newcommand{\Z}{{\mathbb Z}}
\newcommand{\PP}{{\mathcal P}}
\newcommand{\QQ}{{\mathcal Q}}
\newcommand{\EE}{{\mathcal E}}

\newcommand{\Spec}{\operatorname{Spec}}
\newcommand{\D}{{D}}
\newcommand{\R}{\mathbf R}
\renewcommand{\L}{\mathbf L}
\newcommand{\LL}{\mathbf L}
\newcommand{\OO}{{\mathscr O}}
\newcommand{\Ltensor}{\mathbin{\buildrel{\mathbf L}\over{\tensor}}}
\newcommand{\Ext}{\operatorname{Ext}}
\newcommand{\Hilb}{\operatorname{Hilb}}
\newcommand{\Ob}{\operatorname{Ob}}

\begin{document}
\normalsize
\title[]{Equivalences of triangulated categories and Fourier-Mukai transforms}
\author[]{Tom Bridgeland}
\date{\today}

\begin{abstract}
We give a condition for an exact functor
between triangulated categories to be an equivalence. Applications to Fourier-Mukai transforms are discussed. In particular we obtain a
large number of such transforms for K3 surfaces.
\end{abstract}

\maketitle

\section{Introduction}

Let $X$ and $Y$ be smooth projective varieties of the same dimension, and let $\PP$ be a
vector bundle on $X\times Y$. Define a
functor
$$F:D(Y)\lra D(X)$$
between the derived categories of sheaves
on $Y$ and $X$ by the formula
\begin{equation*}
\label{paddle}
F(-)=\R\pi_{X,*}(\PP\tensor\pi_Y^*(-)),
\end{equation*}
where $X\lLa{\pi_X}X\times Y\lRa{\pi_Y} Y$ are the projections maps.
Functors of this type which are equivalences of categories are called
{\it Fourier-Mukai transforms}, and have proved to be powerful tools for studying moduli spaces of vector
bundles [4],[5],[11].

\smallskip

A vector bundle $\PP$ on $X\times Y$ is called {\it strongly simple} over $Y$ if for each
point $y\in Y$, the bundle $\PP_y$ on $X$ is simple, and if for any
two distinct points $y_1,y_2$ of $Y$, and any integer $i$, one has
$$\Ext^{i}_X(\PP_{y_1},\PP_{y_2})=0.$$
One might think of the family $\{\PP_y:y\in Y\}$ as an
`orthonormal' set of bundles on $X$.

\smallskip

The following basic result allows one to construct many examples of Fourier-Mukai
transforms.

\begin{thm}
\label{basic}
The functor $F$ is fully faithful if, and only if, $\PP$ is
strongly simple over $Y$. It is an equivalence of categories precisely
when one also has $\PP_y=\PP_y\tensor\omega_X$ for all $y\in Y$.
\end{thm}

The first statement is well-known [3],[8], but the second part is
new. In this paper we shall prove Theorem \ref{basic}, along with some
more general results concerning exact functors between triangulated categories.

\smallskip

As an example of the use of Therorem \ref{basic}, we have

\begin{cor}
Let $X$ be an algebraic K3 surface and let $Y$ be a fine, compact, 2-dimensional moduli
space of stable vector
bundles on $X$. Then $Y$ is also a K3 surface, and if $\PP$
is a universal bundle on $X\times Y$, the functor $F$ is an equivalence of categories.
\end{cor}

\begin{pf}
The fact that $Y$ is a K3 surface is Theorem 1.4 of [12]. Since
$\omega_X$ is trivial, it is enough to check that $\PP$ is strongly
simple over $Y$. This follows from [12], Proposition 3.12, because any
stable sheaf which moves in a 2-dimensional moduli is semi-rigid.
\qed
\end{pf}

\begin{nota}
All our schemes will be Noetherian schemes. A
sheaf on a scheme $X$ will
mean a coherent $\OO_X$-module, and a point of $X$ will mean a closed
point. If $x$ is a point of $X$ then $\OO_x$ denotes the structure
sheaf of $x$ with reduced scheme structure.

If $a$, $b$ are objects of a triangulated category $\A$, put
$$\Hom_{\A}^i(a,b)=\Hom_{\A}(a,T^ib),$$
where $T:\A\to\A$ is the translation functor.

If $X$ is a scheme, $\D(X)$
will denote the bounded derived category of sheaves on
$X$. For an object $E$ of $\D(X)$, let
$$E^{\dual}=\R\lHom_{\OO_X}(E,\OO_X).$$
We shall write $\HH^i(E)$ for the $i$th
cohomology sheaf of $E$, and $E[n]$ for the object obtained by
shifting $E$ to the left by $n$ places. We say that $E$ is a sheaf if
$\HH^i(E)=0$ when $i\neq 0$.

If $f:X\to Y$ is a morphism of schemes, and $E$ is an object of
$D(Y)$, $\L_pf^*(E)$ denotes the $(-p)$th cohomology object of $\L f^*(E)$.

\end{nota}

\begin{ack}
I learnt a great deal from the papers of Bondal and Orlov. I am also
grateful to my supervisor, Antony Maciocia, for all his help
and encouragement. Finally I would like to thank the referee, who
pointed out a major error in an earlier version of the paper.
\end{ack}


\section{Fully faithful functors}

In this section we give a general criterion for an exact functor
between triangulated categories to be fully faithful. Its proof is
very similar to that of [14], Lemma 2.15.

\begin{defn}
Let $\A$ be a triangulated category. A subclass $\om$ of the objects of
$\A$ will be called a {\it spanning class} for $\A$, if for any object $a$ of
$\A$
\begin{align*}
&\Hom_{\A}^i(\oom,a)=0\quad\forall \oom\in \om\quad\forall
i\in\Z\quad\implies a\isom 0, \\
&\Hom_{\A}^i(a,\oom)=0\quad\forall \oom\in \om\quad\forall
i\in\Z\quad\implies a\isom 0.
\end{align*}
\end{defn}

\begin{eg}
\label{easy}
If $X$ is a smooth projective variety, then
the set 
$$\om=\{\OO_x:x\in X\}$$
is a spanning class for $\A=\D(X)$.
\end{eg}

\begin{pf}
For any object $a$ of $\A=\D(X)$, and any $x\in X$, there is a spectral sequence
$$E^{p,q}_2=\Ext^p_X(\HH^{-q}(a),\OO_x)\implies\Hom_{\A}^{p+q}(a,\OO_x).$$
If $a$ is non-zero, let $q_0$ be the maximal value of $q$ such that $\HH^q(a)$ is non-zero,
and assume that $x$ is a closed point in the support of $\HH^{q_0}(a)$. Then there
is a non-zero element of $E^{0,-q_0}_2$
which survives to give an element of
$\Hom_{\A}^{q_0}(a,\OO_x)$. Serre duality then gives a non-zero
element of $\Hom_{\A}^i(\OO_x,a)$, where $i=\dim X-q_0$.
\qed
\end{pf}

\begin{thm}
\label{one}
Let $\A$ and $\B$ be triangulated categories and let $F:\A\lra\B$ be an exact
functor with a left and a right adjoint. Then $F$ is fully faithful
if, and only if, there exists a spanning class
$\om$ for
$\A$, such that for all elements $\oom_1,\oom_2$ of $\om$, and all integers
$i$, the homomorphism
$$F:\Hom_{\A}^i(\oom_1,\oom_2)\lra\Hom_{\A}^i(F\oom_1,F\oom_2)$$
is an isomorphism.
\end{thm}

\begin{pf}
One implication is clear, so let us assume the existence of $\om$ and
prove that $F$ is fully faithful.

Let $H:\B\lra\A$ be a right adjoint of $F$ and let
$$\eta:1_{\A}\lra H\circ F,\qquad\eps:F\circ H\lra 1_{\B},$$
be the unit and counit respectively of the adjunction
$F\dashv H$. Similarly, let $G:\B\lra\A$ be a left adjoint of $F$ and let
$$\zeta:1_{\B}\lra F\circ G,\qquad\delta:G\circ F\lra 1_{\A},$$
be the unit and counit of $G\dashv F$. Note that by [14], Lemma 1.2,
$G$ and $H$ are also exact functors.

For any pair of objects $a$ and $b$ of $\A$, and any integer $i$, there is a commutative
diagram of group homomorphisms

\begin{equation}
\label{clear}
\begin{CD}
\Hom_{\A}^i(a,b)  @>\eta(b)_* >> \Hom_{\A}^i(a,HFb) \\
@V\delta(a)^* VV                   @VV\beta V \\
\Hom_{\A}^i(GFa,b) @>\alpha >> \Hom_{\B}^i(Fa,Fb)
\end{CD}
\end{equation}
in which $\alpha=\zeta(Fa)^*\circ F$ and
$\beta=\eps(Fb)_*\circ F$ are isomorphisms, and the common diagonal is the map
$$F:\Hom_{\A}^i(a,b)\lra\Hom_{\B}^i(Fa,Fb).$$
When $a$ and $b$ are elements of $\om$ this
map is an isomorphism (by hypothesis), so all the maps
in (\ref{clear}) are isomorphisms.

First we show that for any object $a$ in $\om$, the morphism
$\delta(a)$ is an isomorphism. To see this embed $\delta(a)$ in a
triangle of $\A$:
$$GFa\lRa{\delta(a)} a\lra c\lra T(GFa).$$
For any object $b$ of $\om$ we can apply the functor $\Hom_{\A}(-,b)$
to this triangle and obtain a long exact sequence of groups
\begin{align*}
\cdots\lla\Hom_{\A}(GFa,b)\lLa{\delta(a)^*}\Hom_{\A}(a,b)&\lla\Hom_{\A}(c,b)\lla
\\
\lla\Hom_{\A}^{-1}(GFa,b)&\lLa{\delta(a)^*}\Hom^{-1}_{\A}(a,b)\lla\cdots
\end{align*}
But since the maps $\delta(a)^*$ are all isomorphisms, this implies that
$$\Hom_{\A}^i(c,b)=0\quad\forall b\in\om\quad\forall i\in\Z,$$
so $c\isom 0$ and $\delta(a)$ is an isomorphism.

Now take an object $b$ of $\A$, embed the morphism $\eta(b)$ in a
triangle
$$b\lRa{\eta(b)} HFb\lra c\lra Tb,$$
and apply the functor $\Hom_{\A}(a,-)$ with $a\in\om$. The
homomorphisms
$$\eta(a)_*:\Hom_{\A}^i(a,b)\lra\Hom_{\A}^i(a,HFb)$$
appearing in the resulting long exact sequence are isomorphisms
because of the commuting diagram (\ref{clear}) and the fact that $\delta(a)^*$ is an isomorphism. Arguing as above we conclude that
$c\isom 0$ and hence that $\eta(b)$ is an isomorphism. Since $b$ was
arbitrary, this is enough to show that $F$ is fully faithful.
\qed
\end{pf}


\section{Equivalences of triangulated categories}

Here we give a condition for a fully faithful exact functor between
triangulated categories to be an equivalence. We refer to [9], VIII.2 for the
notion of biproducts in an additive category.

\begin{defn}
A triangulated category $\A$ will be called {\it indecomposable} if
whenever $\A_1$ and $\A_2$ are full subcategories of $\A$ satisfying

(a) for every object $a$ of $\A$ there exist objects $a_j\in\Ob(\A_j)$
such that $a$ is a biproduct of $a_1$ and $a_2$,

(b) for any pair of objects $a_j\in\Ob(\A_j)$,
$$\Hom_{\A}^i(a_1,a_2)=\Hom_{\A}^i(a_2,a_1)=0\quad\forall i\in\Z,$$
then there exists $j$ such that $a\isom 0$ for all $a\in\Ob(\A_j)$.
\end{defn}

\begin{eg}
\label{notdone}
If $X$ is a scheme then $\D(X)$ is
indecomposable if and only if $X$ is connected.
\end{eg}

\begin{pf}
We suppose that $X$ is connected and prove that $\A=\D(X)$ is
indecomposable. The (easy) converse is left to the reader.

Suppose $\A_1$ and $\A_2$ are full subcategories of $\A$
satisfying conditions (a) and (b) of the definition. For any
integral closed subscheme $Y$ of $X$, the sheaf $\OO_Y$ is
indecomposable, and is therefore isomorphic to some object of $\A_j$, $j=1$
or 2. For any point $y\in Y$ we must then have that
$\OO_y$ is also isomorphic to an object of $\A_j$, since otherwise
$(b)$ would imply that $\Hom_{\A}(\OO_Y,\OO_y)=0,$
which is not the case.

Let $X_j$ be the union of those $Y$ such that $\OO_Y$ is
isomorphic to an object of $\A_j$. Then $X_1$ and $X_2$ are closed
subsets of $X$ and $X=X_1\cup X_2$. If a point $x\in X$
lies in $X_1$ and $X_2$ then $\OO_x$
is isomorphic to an object of
$\A_1$ and to an object of $\A_2$. This contradicts (b). Thus the
union is disjoint, and the fact that $X$ is connected implies that one
of the $X_j$ (without loss of generality $X_2$) is empty. But then (b)
implies that for any object $a$ of $\A_2$ one has
$$\Hom_{\A}^i(a,\OO_x)=0\quad\forall i\in\Z\quad\forall x\in X,$$
and hence, by the argument of Example \ref{easy}, $a\isom 0$. This completes the proof.
\qed
\end{pf}

\begin{thm}
\label{two}
Let $\A$ and $\B$ be triangulated categories and let $F:\A\lra\B$ be a
fully faithful exact functor. Suppose that $\B$ is indecomposable, and
that not every object of $\A$ is isomorphic to 0. Then $F$ is an equivalence
of categories if, and only if, $F$ has a left adjoint $G$ and a right
adjoint $H$ such that for any object $b$ of $\B$,
$$Hb\isom 0\implies Gb\isom 0.$$
\end{thm}

\begin{pf}
If $F$ is an equivalence then any quasi-inverse of $F$ is a left and right
adjoint for $F$.
For the converse take an object $b$ of $\B$ and (with notation as in
Theorem \ref{one}) embed the morphism $\eps(b)$ in a
triangle of $\B$:
$$FHb\lRa{\eps(b)} b\lra c\lra T(FHb).$$
Applying $H$ one sees that $Hc\isom 0$, because the fact that $F$ is
fully faithful implies that the morphism $H(\eps(b))$ is an
isomorphism. Define full subcategories $\B_1$ and $\B_2$ of $\B$
consisting of objects satisfying $FHb\isom b$ and $Hb\isom 0$
respectively. Now our hypothesis implies that
$$\Hom_{\B}^i(b_1,b_2)=\Hom_{\B}^i(b_2,b_1)=0\quad\forall i\in\Z,$$
whenever $b_j\in\B_j$. Furthermore, the lemma below applied to the
triangle above shows that every object of $\B$ is a biproduct
$b_1\oplus b_2$. Since $\B$ is indecomposable we must have
$$Hc\isom 0\implies c\isom 0,$$
for any $c\in\Ob(\B)$, so the morphism $\eps(b)$ appearing above is an
isomorphism. Since $b$ was arbitrary,
$F\circ H\isom 1_{\B}$ and $F$ is an equivalence. 
\qed
\end{pf}

\begin{lemma}
Let $\A$ be a triangulated category and let
$$a_1\lRa{i_1} b\lRa{p_2} a_2\lRa{0} Ta_1,$$
be a triangle of $\A$.
Then $b$ is a biproduct of $a_1$ and $a_2$ in $\A$.
\end{lemma}

\begin{pf}
Applying the functors $\Hom_{\A}(-,a_1)$ and $\Hom_{\A}(a_2,-)$, one
obtains morphisms $p_1:b\to a_1$ and $i_2:a_2\to b$, such that
$p_1\circ i_1=1_{a_1}$ and $p_2\circ i_2=1_{a_2}$. The
composition $p_2\circ i_1$ is always 0 and replacing $i_2$ by
$i_2-i_1\circ p_1\circ i_2$, we can assume that
$p_1\circ i_2=0$. Then ([9], VIII.2) it is enough to check that the endomorphism of
$b$ given by
$$\phi=1_b-i_1\circ p_1-i_2\circ p_2$$
is the zero map. But this follows from the fact that
$p_1\circ\phi=p_2\circ\phi=0$.
\qed
\end{pf}


\section{Integral functors}

Throughout this section $X$ and $Y$ are smooth projective varieties
over an algebraically closed field $k$ of characteristic zero, and $\PP$ is an object of
$\D(X\times Y)$. $F$ denotes the exact functor
$$\Phi_{Y\to X}^{\PP}:D(Y)\lra D(X)$$
defined by the formula
$$\Phi_{Y\to X}^{\PP}(-)=\R\pi_{X,*}(\PP\Ltensor\pi_Y^*(-)).$$
Following Mukai, we call $F$ an {\it integral functor}. Here we derive various general
properties of such functors. Most of these appeared in some form in
the original papers of Mukai on Abelian varieties [10],[11].

\smallskip

Given a scheme $S$, one can define a relative version of
$F$ over $S$. This is the functor
$$F_S:\D(S\times Y)\lra\D(S\times
X),$$
given by the formula
$$F_S(-)=\R\pi_{S\times Y,*}(\PP_S\Ltensor\pi_{S\times X}^*(-)),$$
where $S\times X\lLa{\pi_{S\times X}} S\times X\times
Y\lRa{\pi_{S\times Y}} S\times Y$ are the projection maps, and $\PP_S$ is
the pull-back of $\PP$ to $S\times X\times Y$.

\smallskip

The following result is similar to [11], Proposition 1.3.

\begin{lemma}
\label{basechange}
Let $g:T\to S$ be a morphism of schemes, and let $E$ be an object of
$D(S\times Y)$, of finite tor-dimension over $S$. Then there is an
isomorphism
$$F_T\circ\L(g\times 1_Y)^*(E)\isom\L(g\times 1_X)^*\circ F_S(E).$$
\end{lemma}

\begin{pf}
One needs to base-change around the diagram
\begin{equation*}
\begin{CD}
T\times X\times Y @>(g\times 1_{X\times Y}) >> S\times X\times Y \\
@V \pi_{T\times X} VV                   @VV \pi_{S\times X} V \\
T\times X @>(g\times 1_X) >>  S\times X
\end{CD}
\end{equation*}
This is justified by the same argument used to prove [3], Lemma 1.3.
\qed
\end{pf}

We can now show that integral functors preserve families of
sheaves. It is this property which makes them useful for studying
moduli problems. See also [11], Theorem 1.6.

\begin{prop}
Let $S$ be a scheme, and $\EE$ a sheaf on $S\times Y$, flat
over $S$. Suppose that for each $s\in S$, $F(\EE_s)$ is a sheaf on
$X$. Then there is a sheaf $\Hat{\EE}$ on $S\times X$, flat over $S$, such that for every
$s\in S$, $\Hat{\EE}_s=F(\EE_s)$.
\end{prop}

\begin{pf}
Let $\Hat{\EE}=F_S(\EE)$, and take a point $s\in S$. Applying
Lemma \ref{basechange} with $T=\{s\}$, we see
that the derived restriction of $\Hat{\EE}$ to the fibre
$X\times\{s\}$ is just $F(\EE_s)$. The following lemma then shows that
$\Hat{\EE}$ is a sheaf on $S\times X$, flat over $S$.
\qed
\end{pf}

\begin{lemma}
\label{flat}
Let $\pi:S\to T$ be a morphism of schemes,\footnote{There is an error in the published version: we must also assume that $\pi$ is flat. I'm grateful to Chris Seaman for pointing this out.}  and for each
point $t\in T$, let $i_t:S_t\to S$ denote the inclusion of the fibre
$\pi^{-1}(t)$. Let $\EE$
be an object of $D(S)$, such that for all $t\in T$, $\L i_t^*(\EE)$ is
a sheaf on $S_t$. Then $\EE$ is a sheaf on $S$, flat over $T$.
\end{lemma}

\begin{pf}
For each point $t\in T$, consider the
hypercohomology spectral sequence
$$E^{p,q}_2=\L_{-p}i_t^*(\HH^q(\EE))\implies\L_{-(p+q)}i_t^*(\EE).$$
By assumption, the right-hand side is zero unless $p+q=0$.
If $q_0$ is the largest $q$ such that $\HH^q(\EE)\neq 0$, then
$E^{0,q_0}_2$ survives in the spectral sequence for some $t\in T$,
so $q_0=0$. Now $\HH^0(\EE)$ must be flat over $T$ since otherwise\footnote{Here we use the local criterion for flatness: recall that all our schemes are assumed to be Noetherian.}
$E^{-1,0}_2$ would survive for some $t\in T$. Finally, suppose
$\HH^q(\EE)\neq 0$ for some $q<0$. Then we can find a largest such
$q$, and this gives an element of
$E^{0,q}_2$ which survives. Hence $\HH^q(\EE)=0$ unless $q=0$ and
$\EE$ is a sheaf, flat over $T$.
\qed
\end{pf}

In the next section we shall need 

\begin{lemma}
\label{ks}
Suppose that $\PP$ is a sheaf on $X\times Y$, flat over $Y$, and fix a
point $y\in Y$. Then the
homomorphism
\begin{equation}
\label{hee}
F:\Ext^1_Y(\OO_y,\OO_y)\lra\Ext^1_X(\PP_y,\PP_y),
\end{equation}
is the Kodaira-Spencer map for the family $\PP$ at the point $y$, if we identify the
first space with the tangent space to $Y$ at $y$ in the usual way.
\end{lemma}

\begin{pf}
Let $D=\Spec k[\epsilon]/{\epsilon^2}$ denote the double
point. We identify the tangent
space $\T_y Y$ to $Y$ at $y$ with the set of morphisms $D\to
Y$, taking the closed point of $D$ to $y$. Given such a morphism $f$,
we can pull $\PP$ back, and obtain a deformation of the
sheaf $\PP_y$ on $X$, with base $D$. The set of such
deformations is identified with $\Ext^1_X(\PP_y,\PP_y)$, and the
Kodaira-Spencer map is the resulting linear map
$$\T_y Y\lra \Ext^1_X(\PP_y,\PP_y).$$

Returning to our homomorphism (\ref{hee}), note that we can identify the
domain with the set of deformations of $\OO_y$ over $D$, and the image
with the set of deformations of $\PP_y$ over $D$. If we do this, it is
easy to see that the
map $F$ is just given by applying the functor $F_D$.

Given an element $f:D\to Y$ of $\T_y Y$, the corresponding deformation
of $\OO_y$ over $D$ is obtained by pulling-back the
family $\OO_{\Delta}$ on $Y\times Y$ to $D\times Y$ using $f$ (here
$\Delta$ denotes the diagonal in $Y\times Y$). By Lemma \ref{basechange}, if we then apply $F_D$, we get the same
result as if we first applied $F_Y$,
which gives the sheaf $\PP$ on $X\times Y$, and then pulled-back via $f$. But this is
the Kodaira-Spencer map for the family $\PP$.
\qed
\end{pf}

The following result is well-known. Its proof is a straightforward application of Grothendieck-Verdier
duality (see [8], Proposition 3.1 or [3], Lemma 1.2).

\begin{lemma}
\label{adjoint}
The functors
$$G=\Phi_{X\to Y}^{\PP^{\dual}\tensor\pi_X^*\oom_X[\dim X]},\qquad
H=\Phi_{X\to Y}^{\PP^{\dual}\tensor\pi_Y^*\oom_Y[\dim Y]},$$
are
left and right adjoints for $F$ respectively.
\qed
\end{lemma}


\section{Applications to Fourier-Mukai transforms}

As in the last section we fix smooth projective varieties $X$ and $Y$
over an algebraically closed field $k$ of characteristic zero, and an object $\PP$ of $D(X\times Y)$. $F$ denotes the corresponding
functor $\Phi^{\PP}_{Y\to X}$. The following theorem was first proved
by A.I. Bondal and D.O. Orlov, using ideas of Mukai.

\begin{thm}
\label{Borlov}
{\rm ([3])}
The functor $F$ is fully faithful if, and only if, for each point
$y\in Y$,
$$\Hom_{\D(X)}(F\OO_{y},F\OO_{y})=k,$$
and for each pair of points $y_1,y_2\in Y$, and each integer $i$,
$$\Hom_{D(X)}^i(F\OO_{y_1},F\OO_{y_2})=0\mbox{ unless }y_1=y_2\mbox{ and
}0\leq i\leq\dim Y.$$
\end{thm}

\begin{pf}
We must show that for any point $y$ of
$Y$, and any integer $i$, the homomorphism
$$F:\Hom_{\D(Y)}^i(\OO_y,\OO_y)\lra\Hom_{\D(X)}^i(F\OO_y,F\OO_y)$$
is an isomorphism. Theorem \ref{one} will then give the
result.
By the commutative diagram (\ref{clear}) it will be enough to show that
$\delta(\OO_y)$ is an isomorphism. In fact it will be enough to show
that $GF\OO_y\isom \OO_y$, because then $\delta(\OO_y)$ must be either an
isomorphism or zero, and the latter is impossible, because
$F(\delta(\OO_y))$ has a left-inverse.

For any point $z\stackrel{i_z}{\hookrightarrow}Y$, there are isomorphisms of
vector spaces
$$\LL_p
i_z^*(GF\OO_y)\isom\Hom^p_{\D(Y)}(GF\OO_y,\OO_z)\isom\Hom_{\D(X)}^p(F\OO_y,F\OO_z)$$
coming from the adjunctions $i_z^*\dashv i_{z,*}$ ([7], Corollary
5.11), and $G\dashv F$. Thus, by [3], Proposition 1.5, $GF\OO_y$ is a
sheaf supported at the point $y$. Furthermore, there is a unique
morphism $GF\OO_y\to\OO_y$. If $K$ is the kernel of this morphism, one
has a short exact sequence
$$0\lra K\lra GF\OO_y\lRa{\delta(\OO_y)}\OO_y\lra 0,$$
and we must show that $K=0$. Applying the functor
$\Hom_{\D(Y)}(-,\OO_y)$, and using the diagram (\ref{clear}), it will be enough to
show that the homomorphism
\begin{equation}
\label{hope}
F:\Hom_{D(Y)}^1(\OO_y,\OO_y)\lra\Hom^1_{D(X)}(F\OO_y,F\OO_y),
\end{equation}
is injective.

By [10], Proposition 1.3, $GF=\Phi_{Y\to Y}^{\QQ}$ for some object $\QQ$ of
$\D(Y\times Y)$. Since $GF\OO_y$ is a sheaf for all $y\in Y$, Lemma \ref{flat} shows that $\QQ$ is in fact a
sheaf, flat over $Y$. Furthermore, by Lemma \ref{ks}, the map
$$GF:\Hom^1_{D(Y)}(\OO_y,\OO_y)\lra\Hom^1_{D(Y)}(GF\OO_y,GF\OO_y),$$
is given by the
Kodaira-Spencer map for the family $\QQ$ at the point $y$. The following two lemmas
show that this map is injective. Clearly the map
(\ref{hope}) must also be injective.
\qed
\end{pf}

\begin{lemma}
Let $Y$ be a projective variety over $k$, and let $Q$ be a sheaf on $Y$
supported at a point $y\in Y$. Suppose that
$$\Hom_Y(Q,\OO_y)=k.$$
Then $Q$ is the structure sheaf of a zero-dimensional closed subscheme of $Y$.
\end{lemma}

\begin{pf}
There exists a short exact sequence
$$0\lra P\lra Q\lRa{g} \OO_y\lra 0.$$
Suppose $f:\OO_Y\to Q$ is a non-surjective morphism of sheaves.
Considering the cokernel of $f$ shows that there is a non-zero morphism $h:Q\to
\OO_y$ such that $h\circ f=0$. But by hypothesis $h$ must be a
multiple of $g$, so one must have $g\circ f=0$, hence $f$ comes from a
morphism $\OO_Y\to P$. Now
$$\dim_k H^0(Y,P)=\chi(P)<\chi(Q)=\dim_k H^0(Y,Q),$$
so there must be a morphism $\OO_Y\to Q$ which is surjective.
\qed
\end{pf}

\begin{lemma}
Let $S$ and $Y$ be varieties over $k$, with $Y$ projective. Let $\QQ$ be a
sheaf on $S\times Y$, flat over $S$, such that for each $s\in S$,
$\QQ_s$ is the structure sheaf of a zero-dimensional closed subscheme of $Y$. Suppose
also that for all pairs of points $s_1,s_2\in S$
\begin{equation}
\label{emmy}
\QQ_{s_1}\isom \QQ_{s_2}\implies s_1=s_2.
\end{equation}
Then there exists a point $s\in S$, such that the Kodaira-Spencer map for the family $\QQ$ at
$s$ is injective.
\end{lemma}

\begin{pf}
Firstly, we may suppose that $S$ is affine. Fix a point $s\in S$, and let $\pi:S\times Y\to S$
be the projection map. By the theorem on
cohomology and base-change, the
natural map
$$H^0(S\times Y,\QQ)\to H^0(Y,\QQ_y)$$
is surjective, so we can find a section $g:\OO_{S\times Y}\to \QQ$
such that the restriction $g_s:\OO_Y\to \QQ_s$ is surjective. Passing
to an open subset of $S$ we can assume that $g$ is surjective, so that
$\QQ$ is the structure sheaf of a closed subscheme of $S\times Y$.

Let $P$ be the (constant) Hilbert polynomial of the sheaf $\QQ_s$ on
$Y$. By the general existence theorem for Hilbert schemes [6], there is a scheme $\Hilb^P(Y)$ representing the functor which assigns to a
scheme $S$ the set of $S$-flat quotients $\QQ$ of $\OO_{S\times Y}$
with Hilbert polynomial $P$. Let $\EE$
be the universal quotient on $\Hilb^P(Y)\times Y$. Then there is a morphism
$f:S\to \Hilb^P(Y)$ such that $\QQ=(f\times 1_Y)^*(\EE)$. The Kodaira-Spencer
map for the family $\QQ$ at $s\in S$ is obtained by composing the
Kodaira-Spencer map for the family $\EE$ at $f(s)$ with the
differential
$$\T_s(f):\T_s S\lra \T_{f(s)} \Hilb^P(Y).$$
Now condition (\ref{emmy}) implies that the morphism $f$ is injective
on points. Let $S'$ be its scheme-theoretic image. Since we are in
characteristic zero we can assume $S$,$S'$ are non-singular, and
$f:S\to S'$ is smooth. This implies that for some $s\in S$, $\T_s (f)$ is injective. Finally, the fact that the
Kodaira-Spencer map for the family $\EE$ is injective is a consequence
of the universal property of $\EE$. This completes the proof.
\qed
\end{pf}

Theorem \ref{two} allows us to say when $F$ is an equivalence.

\begin{thm}
\label{mine}
Suppose $F$ is fully faithful. Then $F$ is an equivalence if, and only
if, for every point $y\in Y$,
\begin{equation}
\label{last}
F\OO_y\otimes\omega_X\isom F\OO_y.
\end{equation}
\end{thm}

\begin{pf}
Let $G$ and $H$ denote the left and right adjoint functors of $F$
respectively. Suppose first that $F$ is an equivalence. Then $G$ and $H$ are
both quasi-inverses for $F$, so for any $y\in Y$,
$$G(F\OO_y)\isom H(F\OO_y)\isom\OO_y.$$
From the formulas for $G$ and $H$ given in Lemma \ref{adjoint},
$$G(F\OO_y)\isom G(F\OO_y)\otimes\omega_Y \isom H(F\OO_y\otimes\omega_X)[\dim
X-\dim Y].$$
But $G$ is an equivalence, so one
concludes that $X$ and $Y$ have the same dimension, and there is an
isomorphism (\ref{last}).

For the converse, let $X$ have dimension $n$,
and suppose that (\ref{last}) holds for all $y\in Y$. Take an object $b$ of $\D(X)$
such that $Hb\isom 0$. For any point $y\in Y$, and any integer
$i$,
\begin{align*}
\Hom^i_{\D(Y)}(Gb,\OO_y)&=\Hom^i_{\D(X)}(b,F\OO_y)=\Hom^i_{\D(X)}(b,F\OO_y\otimes\omega_X)\\
&=\Hom^{n-i}_{\D(X)}(F\OO_y,b)^{\dual}=\Hom^{n-i}_{\D(Y)}(\OO_y,Hb)^{\dual}=0,
\end{align*}
so by Example \ref{easy}, $Gb\isom 0$. Applying Theorem \ref{two} completes
the proof.
\qed
\end{pf}

Finally note that Theorems \ref{Borlov} and \ref{mine} imply Theorem
\ref{basic} in the special case when $\PP$ is a vector bundle on
$X\times Y$.

\smallskip

\small

\section*{References}

[1] {\it C. Bartocci, U. Bruzzo, D. Hernandez Ruiperez,} A
Fourier-Mukai transform for stable bundles on K3 surfaces, J. reine
angew. Math. {\bf 486} (1997), 1-16.

[2] {\it A.I. Bondal,} Representations of associative algebras and
coherent sheaves, Math. USSR Izv. {\bf 34} (1990), 23-42.

[3] {\it A.I. Bondal, D.O. Orlov,} Semiorthogonal decomposition for
algebraic varieties, Preprint alg-geom/9506012.

[4] {\it T. Bridgeland,} Fourier-Mukai transforms for elliptic
surfaces, Preprint alg-geom/9705002, to appear in J. reine angew. Math.

[5] {\it U. Bruzzo, A. Maciocia,} Hilbert schemes of points on some K3
surfaces and Gieseker stable bundles, Math. Proc. Cam. Phil. Soc. {\bf
120} (1996), 255-261.

[6] {\it A. Grothendieck,} Techniques de construction et th\'{e}or\`{e}mes
d'existence en g\'{e}om\'{e}trie alg\'{e}brique, IV : les sch\'{e}mas de Hilbert, Sem.
Bourbaki {\bf 221} (1960).

[7] {\it R. Hartshorne,} Residues and duality, Lect. Notes Math. {\bf
20}, Springer-Verlag, Heidelberg (1966).

[8] {\it A. Maciocia,} Generalized Fourier-Mukai transforms, J. reine
angew. Math. {\bf 480} (1996), 197-211.

[9] {\it S. MacLane,} Categories for the working mathematician, Grad.
Texts Math.
{\bf 5}, Springer-Verlag, Berlin (1971).

[10] {\it S. Mukai,} Duality between $D(X)$ and $D(\Hat X)$ with its
application to Picard sheaves, Nagoya Math. J. {\bf 81} (1981),
153-175.

[11] {\it S. Mukai,} Fourier functor and its application to the moduli
of bundles on an abelian variety, Adv. Pure Math. {\bf 10} (1987),
515-550.
 
[12] {\it S. Mukai,} On the moduli space of bundles on K3 surfaces I,
in: Vector Bundles on Algebraic Varieties, M.F. Atiyah et al., Oxford
University Press (1987), 341-413.

[13] {\it S. Mukai,} Duality of polarized K3 surfaces, to appear in
Proc. of Algebraic Geometry Euroconference (Warwick 1996).

[14] {\it D.O. Orlov,} Equivalences of derived categories and K3
surfaces, Preprint alg-geom/9606006.

\medskip

Department of Mathematics and Statistics, The University of Edinburgh,
King's Buildings, Mayfield Road, Edinburgh, EH9 3JZ, UK.

email: {\tt tab@maths.ed.ac.uk}

\end{document}